\documentclass{conm-p-l}
\numberwithin{equation}{section}

\def\Cbbd{\mathbb{C}}

\def\Nbbd{\mathbb{N}}
\def\Pbbd{\mathbb{P}}

\def\Acal{\mathcal{A}}

\def\Dcal{\mathcal{D}}

\def\Pcal{\mathcal{P}}

\def\Scal{\mathcal{S}}
\def\abf{\mathbf{a}}
\def\ubf{\mathbf{u}}
\def\vbf{\mathbf{v}}
\def\xbf{\mathbf{x}}
\def\ybf{\mathbf{y}}
\def\zbf{\mathbf{z}}
\def\one{\mathbf{1}}
\def\epsbf{\boldsymbol{\epsilon}}

\def\nubf{\boldsymbol{\nu}}
\def\a{\alpha}
\def\d{\delta}
\def\eps{\epsilon}
\def\D{\Delta}
\def\la{\lambda}

\def\dd{\partial}
\def\Ref#1{(\ref{#1})}

\def\endproof{\hfill\rule{2mm}{2mm}}
\def\?{(?)\marginpar{|?}}

\newtheorem{theo}{Theorem}[section]
\newtheorem{prop}{Proposition}[section]
\newtheorem{lemma}{Lemma}[section]

\def\beq{\begin{equation}}
\def\eeq{\end{equation}}
\def\be{\begin{equation*}}
\def\ee{\end{equation*}}
\begin{document}
\title{Factorization of symmetric polynomials}
\author{Vadim B. Kuznetsov}
  \address{Department of Applied Mathematics,
          University of Leeds,
          Leeds LS2 9JT, UK}
  \email{V.B.Kuznetsov@leeds.ac.uk}
\author{Evgeny K. Sklyanin}
\address{Department of Mathematics, University of York,
York YO10 5DD, UK}
\email{eks2@york.ac.uk }
%
%
\begin{abstract}
We construct linear operators factorizing the three bases of symmetric
polynomials:
monomial symmetric functions $m_\lambda(x)$,
elementary symmetric polynomials $E_\lambda(x)$,
and Schur functions $s_\lambda(x)$, into products of univariate polynomials.
\end{abstract}
\maketitle

\section{Introduction}\label{intro}
\noindent

In our recent paper \cite{KMS} written together with V.V.~Mangazeev
we used the lore
of the quantum integrability to provide a new insight
into the theory of symmetric polynomials.
The main result of \cite{KMS}
is a construction of a {\em separation of variables}
for the Jack polynomials. We have explicitly described
an integral operator $\Scal_n$ depending on the parameter $\a$
which maps any Jack polynomial
$P^{(\a)}_\la(\xbf)$ of $n$ variables $(x_1,\ldots,x_n)\equiv\xbf$
labelled by a partition $\la$ of length $n$
into a product $\prod_{j=1}^n q_\la(z_j)$ of univariate polynomials
$q_\la(z)$. Some other integral operators closely related to $\Scal_n$
were also constructed and studied.

In case of the generic Jack polynomials, the proofs
presented in \cite{KMS} are quite involved, and
the simple underlying philosophy of the separation of variables
is somewhat overshadowed by the abundant
technicalities. The purpose of the present paper is to provide a gentler
pedagogical introduction into the subject of \cite{KMS}.
We revise the results of \cite{KMS} for the three degenerate cases
of the Jack polynomials corresponding to the three particular values
of the parameter $\a$. These are the three well-known
bases of symmetric polynomials:
monomial symmetric functions $m_\lambda(x)$,
elementary symmetric polynomials $E_\lambda(x)$,
and Schur functions $s_\lambda(x)$
corresponding respectively to the values $\infty$, $0$, and $1$
of the parameter $\a$. For these cases the constructions of \cite{KMS}
simplify drastically, and we provide an independent elementary proof for each
of the three cases.

Let $\Cbbd[\xbf]$ be the ring of polynomials in
$x_1$,\ldots,$x_n$, and $\Cbbd[\xbf]^{S_n}$ be the corresponding
subring of {\em symmetric} polynomials. The three bases $m_\lambda(\xbf)$,
$E_\lambda(\xbf)$, and $s_\lambda(\xbf)$ in the ring $\Cbbd[\xbf]^{S_n}$
we are going to work with
are all labelled by {\em partitions}
$\la=(\la_1,\la_2,\dots,\la_n)\in\Nbbd^n$, $\la_1\ge\la_2\ge \ldots\ge\la_n\ge 0$
of length $n$.
The {\em weight} $|\la|$ of a partition $\la$ being defined as
\beq
|\la|=\sum_{i=1}^n\la_i,
\label{jack0}
\eeq
the {\em dominance partial ordering} $\preceq$ for two partitions $\mu$ and $\la$
is defined as follows:
\beq
 \mu\preceq \la \quad \Longleftrightarrow \quad
\Bigl\{ |{\mu}|=|{\la}|\,; \quad
\sum_{j=1}^k \mu_j\leq\sum_{j=1}^k \la_j\,, \quad
k=1,\dots,n-1\Bigr\}\,.
\label{jack5}
\eeq

We shall often use the notation $\la_{i,j}\equiv\la_i-\la_j$,
in particular, $\la_{i,i+1}=\la_i-\la_{i+1}$ assuming $\la_{n+1}\equiv0$,
so that $\la_{n,n+1}=\la_n$.
For the rest of this section $P_\la(\xbf)$ will denote a member of any of the three
polynomial families in question.

In our analysis of the three cases we shall follow the unified plan.
First, we describe a family of $n$ commuting differential operators $\{H_j\}_{j=1}^n$
whose joint eigenfunctions are the polynomials $P_\la(\xbf)$:
\beq
    [H_j,H_k]=0, \qquad H_jP_\la=h_j(\la)P_\la, \qquad j=1,\ldots,n,
\label{eq:hamiltonians}
\eeq
with $h_j(\la)$ being the corresponding eigenvalues. Thus, the above spectral problem
puts the polynomials $P_\la(\xbf)$ into the framework of the theory
of quantum integrable systems.

It is convenient for our purposes to replace the standard normalization
of the polynomials $P_\la(\xbf)$ by the condition of the unit value at the
special point $\one\equiv(1,\ldots,1)$. For any polynomial $P_\la(\xbf)$ thus define
the polynomial $\bar{P}_\la(\xbf)\equiv P_\la(\xbf)/P_\la(\one)$, $\bar{P}_\la(\one)=1$.

We introduce then the polynomial $q_\la(z)\in\Cbbd[z]$ by fixing all the arguments
of $\bar{P}_\la(\xbf)$ except one:
\beq
    q_\la(z)=\bar{P}_\la(z,1,\ldots,1), \qquad q_\la(1)=1,
\label{eq:def_q}
\eeq
and the operator
$Q_z:\Cbbd[\xbf]^{S_n}\mapsto\Cbbd[\xbf]^{S_n}\otimes\Cbbd[z]$
by its eigenvectors and eigenvalues:
\beq
     Q_z\bar{P}_\la=q_\la(z)\bar{P}_\la.
\label{eq:defQ-eigv}
\eeq

The crucial point is that the operator $Q_z$ (which we call $Q$-operator)
admits a simple description in terms independent of the basis $\bar{P}_\la$.
Note that, by definition, the operators $Q_z$ commute:
\beq
    [Q_{z_1},Q_{z_2}]=0, \qquad \forall z_1,z_2\in\Cbbd.
\label{eq:commutQ}
\eeq

We show that the polynomials $q_\la(z)$ satisfy an ordinary differential equation
in $z$ with the coefficients depending linearly on all eigenvalues
$h_j(\la)$ of the operators $H_j$. The aforementioned
differential equation in $z$ bears the name {\it separation equation}. In other
words, the eigenvalue $q_\la(z)$ of the $Q$-operator solves the
multiparameter spectral
problem arising in the process of separation of variables.
Therefore, these univariate polynomials become {\it
separated polynomials} in the unified approach to the {\it method of
separation of variables} which is explained below.

The next step is to take the composition of $n$ copies of the $Q$-operator
$Q_{z_1}\ldots$ $Q_{z_n}$ and the linear functional
\beq
\rho_0:\Cbbd[\xbf]^{S_n}\mapsto\Cbbd:f(\xbf)\mapsto f(\one).
\label{eq:def_rho0}
\eeq

As a consequence of the commutativity \Ref{eq:commutQ} of $Q_z$
the resulting operator
\beq
    \Scal_n=\rho_0 Q_{z_1}\ldots Q_{z_n}
\label{eq:def_S}
\eeq
acts from $\Cbbd[\xbf]^{S_n}$ into $\Cbbd[\zbf]^{S_n}$
and factorizes the polynomials $\bar{P}_\la(\xbf)$
\beq
    \Scal_n:\bar{P}_\la(\xbf)\mapsto \prod_{j=1}^n q_\la(z_j)
\label{eq:actSP}
\eeq
into products of the univariate polynomials $q_\la(z)$. In this way,
the original spectral problem \Ref{eq:hamiltonians} gets mapped
into a {\it multiparameter spectral problem} for the factorized
eigenfunctions $\prod_{j=1}^n q_\la(z_j)$.

The expression \Ref{eq:def_S} can be simplified due to the
existence of additional relations between $\rho_0$ and $Q_z$.
Define $\rho_k$, $k=0,\ldots,n-1$ as
\beq
    \rho_k:\Cbbd[x_{k+1},\ldots,x_n]^{S_{n-k}}\mapsto \Cbbd
    :f(x_{k+1},\ldots,x_n)\mapsto f(1,\ldots,1)
\label{eq:def_rhok}
\eeq
and let $\rho_n=1$.
We prove then that there exist operators
\beq
\Acal_k:\Cbbd[x_1,\ldots,x_k]^{S_k}\mapsto
  \Cbbd[x_1,\ldots,x_{k-1}]^{S_{k-1}}\otimes\Cbbd[z_k], \qquad
k=1,\ldots,n
\label{eq:Ak_spaces}
\eeq
such that
\beq
    \rho_{k-1} Q_{z_{k}}=\Acal_{k}\rho_{k}.
\label{eq:rhoQ=Arho}
\eeq
(see \cite{KMS}, Proposition 6.1).

Using the commutation relations \Ref{eq:rhoQ=Arho} we can transform
the original expression  $\rho_0 Q_{z_1}\ldots Q_{z_n}$ for $\Scal_n$
into $\Acal_1\rho_1 Q_{z_2}\ldots Q_{z_n}$, then into
$\Acal_1\Acal_2\rho_2Q_{z_3}\ldots Q_{z_n}$ etc.,
and finally, using $\rho_n\equiv1$, into
\beq
    \Scal_n=\Acal_1\Acal_2\ldots\Acal_n\,,
\label{eq:Achain}
\eeq
see \cite{KMS}, Theorem 6.1.
The factorization \Ref{eq:Achain} of the operator $\Scal_n$
exhibits a triangularity property in the sense that each operator $\Acal_k$
acts on $k$ variables only.

Applying the relations \Ref{eq:rhoQ=Arho} to the polynomial $\bar{P}_\la(\xbf)$
we get a set of remarkable equalities
\beq
    \Acal_k:\bar{P}_\la(x_1,\ldots,x_k,1,\ldots1)\mapsto
    \bar{P}_\la(x_1,\ldots,x_{k-1},1,\ldots,1)q_\la(z_k).
\label{eq:actAP}
\eeq

Note that the equalities \Ref{eq:actAP} \underline{cannot} be used to
define the operators $\Acal_k$ directly because the restricted polynomials
$\bar{P}_\la(x_1,\ldots,x_k,1,\ldots1)$ are not linearly independent in
the space $\Cbbd[x_1,\ldots,x_k]^{S_k}$ and, consequently, \Ref{eq:actAP}
contain a lot of nontrivial identities.

Another application of the operator $Q_z$ is to produce the operators
raising the number of variables in the polynomials $\bar{P}_\la$.
As shown in \cite{KMS} in case of generic Jack polynomials the operator
$Q_z$ taken at $z=0$ can be decomposed uniquely
\beq
    Q_0=Q_0^\prime\Pcal
\label{eq:factQ0}
\eeq
into the projector
\begin{align}\label{eq:Pcal}
   \Pcal&:\Cbbd[x_1,\ldots,x_n]^{S_n}\mapsto\Cbbd[x_1,\ldots,x_{n-1}]^{S_{n-1}} \\
      &:f(x_1,\ldots,x_n)\mapsto f(x_1,\ldots,x_{n-1},0) \notag
\end{align}
and an operator
\beq
    Q_0^\prime:\Cbbd[x_1,\ldots,x_{n-1}]^{S_{n-1}}\mapsto
        \Cbbd[x_1,\ldots,x_n]^{S_n}.
\label{eq:Q0prime}
\eeq

The operator $Q_0^\prime$ can be used then to add an extra variable
$x_n$ to the polynomials $\bar{P}_{\la_1\ldots\la_{n-1}}$
\beq
    Q_0^\prime: \bar{P}_{\la_1\ldots\la_{n-1}}(x_1,\ldots,x_{n-1})
\mapsto \bar{P}_{\la_1\ldots\la_{n-1}0}(x_1,\ldots,x_n).
\eeq

In the next three sections the program sketched above
is fulfilled
for each of the three symmetric polynomial bases:
$m_\la(\xbf)$, $E_\la(\xbf)$, and $s_\la(\xbf)$.
In each case our presentation follows the same routine:
\begin{enumerate}
\item\label{lst:basis}
Basis $P_\la$ and the normalized polynomials $\bar{P}_\la$.
\item\label{lst:hamiltonians}
Commuting differential operators $\{H_j\}_{j=1}^n$ with the eigenfunctions $\bar{P}_\la(\xbf)$.
\item\label{lst:Qop}
$Q$-operator $Q_z$ and its eigenvalues $q_\la(z)$.
Differential equation for $q_\la(z)$.
\item\label{lst:Sop}
Separating operator $\Scal_n$ and its
factorization into the chain of $\Acal$-operators.
\item\label{lst:raise}
Lifting operator $Q_0^\prime$.
\end{enumerate}

\section{Factorizing the basis $m_\la(\xbf)$}\label{factm}
\noindent
\subsection{Basis $m_\la(\xbf)$}

Let $\xbf^{\abf}\equiv x_1^{a_1}\ldots x_n^{a_n}$ for any $\abf=(a_1,\ldots,a_n)$.
The {\em monomial symmetric functions} $m_\la(\xbf)$ are
defined as the sum over all \underline{distinct} permutations $\nubf=(\nu_1,\ldots,\nu_n)$
of the partition $\la=(\la_1,\ldots,\la_n)$
\beq
 m_\la(\xbf)=\sum x_1^{\nu_1}\cdots x_n^{\nu_n}
      =x_1^{\la_1}x_2^{\la_2}\ldots x_n^{\la_n}+\text{permuted terms}
\eeq
and are known to form a basis in $\Cbbd[\xbf]^{S_N}$.
We normalize them introducing
$\bar{m}_\la(\xbf)\equiv m_\la(\xbf)/m_\la(\one)$, see Section \ref{intro}.
Note that $\bar{m}_\la(\xbf)$ can be represented as a sum
\beq
    \bar{m}_\la(\xbf)= \frac{1}{n!}\sum_{\sigma\in S_n} \xbf^{\sigma(\la)},
    \qquad \bar{m}_\la(\one)=1,
\label{eq:mbar}
\eeq
over \underline{all} permutations $\sigma\in S_n$.

The particular functions $m_\la$ for
\be
   \la= \underbrace{1\ldots1}_r\underbrace{0\ldots0}_{n-r}\equiv(1^r0^{n-r}),
   \qquad r=0,\ldots,n,
\ee
are known as {\em elementary symmetric functions}. Equivalently,
$e_r(\xbf)$ is the sum of all products of $r$ distinct variables
$x_i$, so that $e_0(\xbf)=1$ and
\beq
e_r(\xbf)=m_{(1^r0^{n-r})}(\xbf)
   =\sum_{1\leq i_1<\ldots<i_r\leq n}x_{i_1}x_{i_2}\ldots x_{i_r}\,.
\label{eq:def_er}
\eeq
Their generating function is
\beq
w_n(t)\equiv\prod_{i=1}^n(1+tx_i)=\sum_{j=0}^n e_j(\xbf) t^j.
\label{eq:genfun_e}
\eeq

\subsection{Operators $H_j$}
The simple monomials $\xbf^\abf$ form a basis in $\Cbbd[\xbf]$ and are already
factorized. They are joint eigenfunctions of $n$ commuting differential operators
\beq
    D_j=x_j\frac{\dd}{\dd x_j}, \qquad
    D_j \xbf^\abf=a_j \xbf^\abf, \qquad j=1,\ldots,n.
\eeq

The symmetric monomial functions $\bar{m}_\la$, however, are not factorized, and the
corresponding differential operators $H_j$ are obtained as elementary
symmetric functions in $D_j$:
\beq
H_j=e_j(D_1,\ldots,D_n),\qquad
H_j\bar{m}_\la=e_j(\la)\bar{m}_\la,
\label{eq:Hops_m}
\eeq
\be
[H_j,H_k]=0,\qquad \forall j,k=1,\ldots,n.
\ee

\subsection{$Q$-operator}
The $Q$-operator $Q_z$ on $\Cbbd[\xbf]^{S_n}$
is defined in terms of its eigenvectors $\bar{m}_\la$
and eigenvalues $q_\la(z)$
\beq
    Q_z\bar{m}_\la=q_\la(z)\bar{m}_\la, \qquad
    q_\la(z)=\bar{m}_\la(z,1,\ldots,1)=\frac{1}{n}\sum_{j=1}^n z^{\la_j}.
\label{eq:eigvalQ_m}
\eeq

\begin{prop}
The operator $Q_z$ admits an alternative description independent of the basis $\bar{m}_\la$:
\begin{subequations}\label{eq:Qop_m_both}
\begin{align}
   Q_z&=\frac{1}{n}\sum_{j=1}^n \;z^{D_j}, \\
   (Q_zf)(\xbf)&=\frac{1}{n}\sum_{j=1}^n f(x_1,\ldots,x_{j-1},zx_j,x_{j+1},\ldots,x_n).
\label{eq:Qop_m}
\end{align}
\end{subequations}
\end{prop}
{\bf Proof.} The expression \Ref{eq:Qop_m_both} is well defined on any polynomial
$f$, not necessarily symmetric. In particular, for the monomials $\xbf^{\abf}$
we have
\be
   Q_z \xbf^{\la}= \frac{1}{n}\sum_{j=1}^n \;z^{D_j}\xbf^\la
   = \left(\frac{1}{n}\sum_{j=1}^n z^{\la_j}\right)\xbf^{\la}
   =q_\la(z)\xbf^{\la}.
\ee

Averaging the monomials $\xbf^{\la}$ over the permutation group with \Ref{eq:mbar}
we recover \Ref{eq:eigvalQ_m}.
\endproof

\subsection{Separation of variables}
\begin{prop}
For the separating operator $\mathcal{S}_n$ defined
by the formula \Ref{eq:def_S}
and satisfying \Ref{eq:actSP}
for $\bar{P}_\la\equiv \bar{m}_\la$
there exist uniquely defined operators $\Acal_k$ satisfying
the relations \Ref{eq:rhoQ=Arho}.
\end{prop}

{\bf Proof.} Applying the both sides of \Ref{eq:rhoQ=Arho}
to a polynomial $f\in\Cbbd[\xbf]^{S_n}$ and using
the formula \Ref{eq:Qop_m} for $Q_zf$ and the definition \Ref{eq:def_rhok}
of $\rho_k$ we obtain:
\begin{align}\label{eq:rhoQf}
   [\rho_{k-1}Q_zf](x_1,\ldots,x_n)
   =&\frac{1}{n} \,f(zx_1,x_2,\ldots,x_{k-1},1,\ldots,1) \\
    &+\ldots \notag \\
    &+\frac{1}{n} \,f(x_1,\ldots,x_{k-2},zx_{k-1},1,\ldots,1) \notag\\
    &+\frac{1}{n} \,f(x_1,\ldots,x_{k-1},z,1,\ldots,1) \notag\\
    &+\ldots \notag\\
    &+\frac{1}{n} \,f(x_1,\ldots,x_{k-1},1,\ldots,1,z). \notag
\end{align}

The last $(n-k+1)$ terms are all equal due to the symmetry of $f$.
The right-hand side of \Ref{eq:rhoQf} has obviously the form
$\Acal_k\rho_k$ where the operator $\Acal_k$ is described
in terms of the projector operators $\Pbbd_{jk}$,
$j<k$, defined as
\beq
    (\Pbbd_{jk}f)(x_1,\ldots,x_n)
    =f(\ldots,\stackrel{j\mathrm{th}}{x_jx_k},\ldots,\stackrel{k\mathrm{th}}1,\ldots).
\label{eq:defPproj}
\eeq

It is convenient to identify the variable $z_k$ in the target space
\Ref{eq:Ak_spaces} for $\Acal_k$ with $x_k$.
The formula for $\Acal_k$ in terms of $\Pbbd_{jk}$ is then
\beq
\Acal_k=\frac{1}{n}\left(n-k+1+\sum_{j=1}^{k-1}\Pbbd_{jk}\right),\qquad \Acal_1\equiv 1.
\label{eq:Ak_proj}
\eeq
\endproof

Note that the formula \Ref{eq:Ak_proj} allows to extend the operator
$\Acal_k$ to all nonsymmetric polynomials from $\Cbbd[\xbf]$.
It is not obvious at all that product of $\Acal_k$ maps
symmetric polynomials into symmetric. Nevertheless, the above analysis
shows that it is true!

Using the relations
\beq
    \Pbbd_{jk}\Pbbd_{lk}=\Pbbd_{jk},\quad \forall j,l=1,\ldots,k-1,
\eeq
for the projectors $\Pbbd_{jk}$ it is easy to invert $\Acal_k$:
\beq
\Acal_k^{-1}=\frac{1}{n-k+1}\left(n-\sum_{j=1}^{k-1}\Pbbd_{jk}\right),
\qquad k=1,\ldots,n.
\eeq

The above inversion formula for $\Acal_k$ seems to be a peculiarity of
the monomial functions case. Its analogs for the generic Jack
polynomials are unknown.

\subsection{Lifting operator}

Letting $z=0$ in \Ref{eq:Qop_m} and using the symmetry of a polynomial
$f\in\Cbbd[\xbf]^{S_n}$ we can shuffle the arguments $\xbf$ to obtain
\beq
   Q_0: f(\xbf)\mapsto \frac{1}{n} \sum_{j=1}^n \left.f(\xbf)\right|_{x_j=0}
     = \frac{1}{n} \sum_{j=1}^n f(x_1,\ldots,\hat{x}_j,\ldots,x_n,0),
\label{eq:Q0_m}
\eeq
the hat over ${x}_j$ marking the omitted argument. Using the formula \Ref{eq:Q0_m} for
$Q_0$ and \Ref{eq:factQ0} we get the expression for the action of the
operator $Q_0^\prime$ on a polynomial $f\in\Cbbd[x_1,\ldots,x_{n-1}]^{S_{n-1}}$:
\beq
    Q_0^\prime: f(x_1,\ldots,x_{n-1}) \mapsto
    \frac{1}{n} \sum_{j=1}^n f(x_1,\ldots,\hat{x}_j,\ldots,x_n).
\label{eq:Q0prime_m}
\eeq

\begin{prop}
The operator $Q_0^\prime$ defined by \Ref{eq:Q0prime_m} acts on the
basis $\bar{m}$ as follows:
\beq
    Q_0^\prime: \bar{m}_{\la_1\ldots\la_{n-1}}(x_1,\ldots,x_{n-1})
\mapsto \bar{m}_{\la_1\ldots\la_{n-1}0}(x_1,\ldots,x_{n}).
\label{eq:Q0prime_act_m}
\eeq
\end{prop}
{\bf Proof.} The formula \Ref{eq:Q0prime_m} defines $Q_0^\prime$ for
nonsymmetric polynomials from $\Cbbd[x_1,\ldots,x_{n-1}]$ as well.
On the monomials we  have
\be
   Q_0^\prime:x_1^{\la_1}\ldots x_{n-1}^{\la_{n-1}}
\mapsto \frac{1}{n}\sum_{j=1}^n
x_1^{\la_1}\ldots x_{j-1}^{\la_{j-1}}x_j^0 x_{j+1}^{\la_j}\ldots x_n^{\la_{n-1}}.
\ee

After having symmetrized the left-hand side over $S_{n-1}$ we obtain in the
right-hand side the average over $S_n$, hence \Ref{eq:Q0prime_act_m}.
\endproof

\section{Factorizing the basis $E_\la(\xbf)$}\label{factE}
\nopagebreak[4]\subsection{Basis $E_\la(\xbf)$}

For each partition $\la=(\la_1,\dots,\la_n)$ define
the polynomials $E_{\la}(\xbf)$ as
\beq
E_{\la}(\xbf)
  =e_1^{\la_1-\la_2}(\xbf)e_2^{\la_2-\la_3}(\xbf)
\cdots e_n^{\la_n}(\xbf)
\label{eq:def_E}
\eeq
where $e_j$, $j=1,\ldots,n$, are elementary symmetric functions defined in
\Ref{eq:def_er}--\Ref{eq:genfun_e}.

The polynomials $E_\la$ form a basis in $\Cbbd[\xbf]^{S_n}$, see \cite{Macd}.
Noticing that $e_j(\one)=\binom{n}{j}$ we introduce the normalized polynomials
\beq
    \bar{e}_j(\xbf)\equiv \frac{e_j(\xbf)}{e_j(\one)}, \qquad
    \bar{E}_\la(\xbf)\equiv \prod_{j=1}^n \bar{e}_j(\xbf)^{\la_{j,j+1}},
\label{eq:normalE}
\eeq
such that $\bar{e}_j(\one)=1$ and $\bar{E}_\la(\one)=1$.

\subsection{Operators $H_j$}
Note that the polynomials $E_\la(\xbf)$ are already given in a factorized form
by the definition \Ref{eq:def_E}. The corresponding separation of variables
is obtained by taking $\eps_j\equiv e_j(\xbf)$ as independent variables
and using the isomorphism $\Cbbd[\xbf]^{S_n}\simeq\Cbbd[\epsbf]$.
The substitution $e_j(\xbf)=\eps_j$ produces the separating linear operator
\beq
    \Scal_n^{(0)}:\Cbbd[\xbf]^{S_n}\mapsto\Cbbd[\epsbf]
    :E_\la(\xbf)\mapsto \prod_{j=1}^n \eps_j^{\la_{j,j+1}}.
\eeq

The monomials $\Scal_n^{(0)}E_\la$ are obviously the eigenfunctions
of the commuting differential operators $\Dcal_j$ in variables $\epsbf$:
\beq
    \Dcal_j=\eps_j\frac{\dd}{\dd\eps_j}, \qquad
    \Dcal_j\Scal_n^{(0)}E_\la= \la_{j,j+1} \Scal_n^{(0)}E_\la,
    \qquad j=1,\ldots,n.
\eeq

It is an easy exercise to express the operators $\Dcal_j$ in terms of the original
variables $\xbf$. Indeed, using the generating function \Ref{eq:genfun_e}
for $\eps_j$ and
taking the logarithmic differential of the both sides one obtains
\beq
\sum_{i=1}^n \frac{t}{1+tx_i}\,dx_i
=\sum_{j=1}^n \frac{t^j}{\prod_{m=1}^n(1+tx_m)}
\,d\eps_j.
\eeq

Comparing the residues at the poles one gets
\beq
dx_i=\sum_{j=1}^n \frac{(-x_j)^{n-j}}{\prod_{m\neq i}^n(x_m-x_i)} \,d\eps_j,\qquad
i=1,\ldots,n,
\eeq
then
\beq
\frac{\partial}{\partial \eps_j}
=\sum_{i=1}^n \frac{(-x_i)^{n-j}}{\prod_{m\neq i}^n(x_m-x_i)} \,
\frac{\partial}{\partial x_i}\,,\qquad j=1,\ldots,n,
\eeq
and, finally,
\beq
H_j\equiv \bigl(\Scal_n^{(0)}\bigr)^{-1}\Dcal_j\Scal_n^{(0)}
=e_j(\xbf)\sum_{i=1}^n \frac{(-x_i)^{n-j}}{\prod_{m\neq i}^n(x_m-x_i)} \,
\frac{\partial}{\partial x_i}\,,\qquad j=1,\ldots,n.
\label{dkey}
\eeq

As a consequence, the operators $H_j$
mutually commute and have the functions $E_\la(\xbf)$ as their eigenfunctions:
\beq
[H_j,H_k]=0,\qquad
H_j E_\la(\xbf)=\la_{j,j+1} E_\la(\xbf),\qquad j=1,\ldots,n.
\label{spe}
\eeq

The factorizing operator $\Scal_n$ obtained from the constructions of \cite{KMS}
in the limit $\a\rightarrow0$ and described below
is, however, different from $\Scal_n^{(0)}$ described above,
illustrating thus the point that a separation of variables is not unique.
Both operators $\Scal_n$ and $\Scal_n^{(0)}$ correspond to the same set of
commuting differential operators $H_j$ but the separated functions and the
separated equations are different.

\subsection{$Q$-operator}
Define the operator $Q_z$ by its eigenvectors $\bar{E}_\la(\xbf)$ and eigenvalues
\beq
Q_z: \bar{E}_\la\mapsto q_\la(z) \bar{E}_\la,
\eeq
\begin{align}
q_\la(z)&\equiv\bar{E}_\la(z,1,\ldots,1)
        =\prod_{j=1}^n \bar{e}_j(z,1,\ldots,1)^{\la_{j,j+1}} \\
        &=\prod_{j=1}^n \left(1+\frac{z-1}{n}j \right)^{\la_{j,j+1}}. \notag
\end{align}

The polynomial $q_\la(z)$ is the unique polynomial solution of the
differential equation
\beq
    \frac{dq_\la(z)}{dz}=\sum_{j=1}^n\frac{\la_{j,j+1}}{z+\frac{n-j}{j}}q_\la(z)
\eeq
satisfying the normalization condition $q_\la(1)=1$.

The factorized structure of the polynomials $E_\la(\xbf)$
and $q_\la(z)$ implies that the action of the operator $Q_z$ on
polynomials from $\Cbbd[\xbf]^{S_n}$ is determined uniquely by its
action on the elementary symmetric functions $\eps_j$:
\beq
    Q_z:\eps_j\equiv e_j(\xbf)\mapsto \left(1+\frac{z-1}{n}j\right)\eps_j.
\label{eq:actQe}
\eeq

The formula \Ref{eq:actQe}
can be written in terms of the generating function \Ref{eq:genfun_e}
as follows:
\beq
    Q_z: w_n(t)\mapsto \left(1+\frac{z-1}{n}\,t\frac{d}{d t}\right)w_n(t).
\label{eq:Qgenfun}
\eeq

Using the isomorphism $\Cbbd[\xbf]^{S_n}\simeq\Cbbd[\epsbf]$ we observe
that $Q_z$ acts on a polynomial $f\in\Cbbd[\epsbf]$ by substitution
\beq
    Q_z:f(\eps_1,\ldots,\eps_n)\mapsto f\bigl(Q_z(\eps_1),\ldots,Q_z(\eps_n)\bigr).
\eeq

In fact, all operators we are going to construct in this section:
$\Scal_n$, $\Acal_k$, $Q_0^\prime$, share with $Q_z$ the same property
of acting on $\Cbbd[\epsbf]$ by substitutions.

\subsection{Separation of variables}
The separating operator $\Scal_n$ is defined by the formula \Ref{eq:def_S}.
By virtue of \Ref{eq:actSP},
it acts on the basis $\bar{E}_\la(\xbf)$ as follows:
\beq
\Scal_n: \bar{E}_\la(\xbf)\mapsto \prod_{j=1}^n q_\la(z_j).
\label{eq:actSE}
\eeq

As explained above, the formula \Ref{eq:actSE} implies that
the operator $\Scal_n$ acts on $\Cbbd[\epsbf]$ by substitutions:
\beq
    \Scal_n:\eps_j\mapsto \binom{n}{j}
\prod_{i=1}^n \left(1+\frac{z_i-1}{n}j\right)
\eeq
(note the normalization factor $\binom{n}{j}=e_j(\one)$ in the
right-hand side).

\begin{theo}
For the separating operator $\mathcal{S}_n$ defined
above there exist uniquely defined operators $\Acal_k$ satisfying
the relations \Ref{eq:rhoQ=Arho} and thus \Ref{eq:Achain}.
\end{theo}
{\bf Proof.}
Applying the hypothetic operator equality \Ref{eq:rhoQ=Arho}
to the basis $\bar{E}_\la(\xbf)$ we get the system of relations
\Ref{eq:actAP}
for $\bar{P}_\la=\bar{E}_\la$. The existence
and uniqueness of the operators $\Acal_k$ will be established
when we prove that the overdetermined system \Ref{eq:actAP}
has a unique solution. Because of the factorized nature
of the polynomials $\bar{E}_\la(\xbf)$ and $q_\la(\xbf)$
it is sufficient to consider the relations \Ref{eq:actAP}
for the elementary
symmetric functions $e_j(\xbf)$ only. The resulting equations
determining the action of $\Acal_k$ on $\rho_k e_j$ are:
\beq
    \Acal_k\rho_k e_j=\left(1+\frac{z_k-1}{n}j\right)\rho_{k-1}e_j,
  \qquad j=1,\ldots,n.
\label{eq:Akrhokej}
\eeq

The system \Ref{eq:Akrhokej} is apparently overdetermined because it counts $n$
equations for only $k$ independent quantities contained in $\rho_k e_j$
(note that $\rho_k e_j$ depend only on the variables $x_1,\ldots,x_k$).
To prove that the equations \Ref{eq:Akrhokej} are nevertheless consistent
let us sum them with the factor $t^j$ and use the
correspondence \Ref{eq:actQe}--\Ref{eq:Qgenfun} to get
\beq
    \Acal_k \rho_k w_n(t)=
    \left(1+\frac{z_k-1}{n}\,t\frac{d}{d t}\right)\rho_{k-1}w_n(t).
\label{eq:Akrhokwn}
\eeq

{}From the definition \Ref{eq:def_rhok} it follows that
\beq
    \rho_k w_n(t)=(1+t)^{n-k}\prod_{i=1}^k (1+tx_i)
      =(1+t)^{n-k}w_k(t).
\label{eq:rhok_wn}
\eeq

Substituting $\rho_k w_n(t)$ and $\rho_{k-1} w_n(t)$ from
\Ref{eq:rhok_wn} to \Ref{eq:Akrhokwn} and dividing by
$(1+t)^{n-k}$ we get
\beq
   \Acal_k w_k(t)=
   (1+t)^{k-n}\left(1+\frac{z_k-1}{n}\,t\frac{d}{d t}\right)
   (1+t)^{n-k+1}w_{k-1}(t).
\eeq

Simplifying the resulting expression and expanding it in $t$ we finally
get
\begin{subequations}
\begin{align}
    \Acal_k:e_j^{(k)}&\mapsto e_j^{(k-1)}
    \left(1+\frac{z_k-1}{n}j\right) \\
    &+e_{j-1}^{(k-1)}\left(1+\frac{z_k-1}{n}(n-k+j)\right),
    \qquad j=1,\ldots,k-1, \notag \\
   \Acal_k:e_k^{(k)}&\mapsto z_k e_{k-1}^{(k-1)},
\end{align}
\end{subequations}
where $e_j^{(k)}$ is the $j^{\text{th}}$ elementary symmetric function
in variables $x_1,\ldots,x_k$.

Remarkably, because of the cancellation of the factor $(1+t)^{n-k}$,
the resulting system contains only $k$ equations determining thus
uniquely the action of $\Acal_k$ on $e_j^{(k)}$.
\endproof

\subsection{Lifting operator}
\begin{prop}
The operator $Q_z$ taken at $z=0$ admits a unique factorization
\Ref{eq:factQ0} into the projector $\Pcal$ \Ref{eq:Pcal}
and the operator $Q_0^\prime$ \Ref{eq:Q0prime}:
\beq
Q_0^\prime: \bar{E}_{\la_1\ldots\la_{n-1}}(x_1,\ldots,x_{n-1})
\mapsto \bar{E}_{\la_1\ldots\la_{n-1}0}(x_1,\ldots,x_{n}).
\label{eq:Q0prime_E_lift}
\eeq
\end{prop}
{\bf Proof.}
As usual, it is sufficient to consider the action of the studied operators on
the elementary symmetric functions $e_j$ only.
Setting $z=0$ in \Ref{eq:actQe} we get
\be
Q_0:e_j^{(n)}\mapsto \frac{n-j}{n}\,e_j^{(n)}.
\ee

{}From the definition \Ref{eq:Pcal} of the projector $\Pcal$
it follows that
\begin{align*}
    \Pcal&:e_j^{(n)}\mapsto e_j^{(n-1)}, \qquad j=1,\ldots,n-1, \\
    \Pcal&:e_n^{(n)}\mapsto 0.
\end{align*}

The unique factorization $Q_0=Q_0^\prime\Pcal$ is then obvious, with
\be
    Q_0^\prime:e_j^{(n-1)}\mapsto \frac{n-j}{n}\,e_j^{(n)},
    \qquad j=1,\ldots,n-1.
\ee

For the normalized polynomials \Ref{eq:normalE} we get
\be
    Q_0^\prime:\bar{e}_j^{(n-1)}\mapsto \bar{e}_j^{(n)},
    \qquad j=1,\ldots,n-1,
\ee
hence \Ref{eq:Q0prime_E_lift}.
\endproof

\section{Factorizing the basis $s_\lambda(\xbf)$}\label{schur}
\noindent
\subsection{Basis $s_\la(\xbf)$}
Let $\d\equiv(n-1,n-2,\ldots,0)$ and
\beq
     \mu\equiv\la+\d=(\la_1+n-1,\ldots,\la_{n-1}+1,\la_n)
\eeq
for any partition $\la$.
The {\em Schur function} $s_\la(\xbf)$ is defined then as the ratio
of two antisymmetric polynomials
\beq
   s_\la(\xbf)=\frac{\det\{x_i^{\mu_j}\}}{\det\{x_i^{\d_j}\}}
   =\frac{a_\mu(\xbf)}{a_\d(\xbf)}\,,
\label{eq:def_s}
\eeq
where
\beq
    a_\mu(\xbf)= \det\{x_i^{\mu_j}\}
    =\begin{vmatrix}
      x_1^{\mu_1} & x_1^{\mu_2} & \ldots & x_1^{\mu_n} \\
      x_2^{\mu_1} & x_2^{\mu_2} & \ldots & x_2^{\mu_n} \\
             \hdotsfor{4} \\
      x_n^{\mu_1} & x_n^{\mu_2} & \ldots & x_n^{\mu_n}
    \end{vmatrix},
\label{eq:amu}
\eeq
and $a_\d(\xbf)$ is the Vandermonde determinant $\D_n(\xbf)$
\beq
  a_\d(\xbf)
  =\det\{x_i^{\d_j}\}
  =\D_n(\xbf)=\prod_{i<j}(x_i-x_j)
\label{eq:vandermonde}
\eeq
(note that the sign of $\D_n(\xbf)$ differs from \cite{KMS} by the
factor $(-1)^{n(n-1)/2}$). In the formulas like \Ref{eq:def_s}
it is always assumed that in the expression under the
$\det$ sign the index $i$ is the number of the matrix row,
and $j$ of the column.

\begin{lemma}
The Schur function evaluated at $x_{k}=\ldots=x_n=1$
is given by the formula
\begin{subequations}\label{eq:restr_schur}
\beq
    s_\la(x_1,\ldots,x_{k-1},1,\ldots,1)
   =\frac{a_\mu^{(k)}(x_1,\ldots,x_{k-1})}{a_\d^{(k)}(x_1,\ldots,x_{k-1})},
\qquad k=1,\ldots,n,
\label{eq:s=a/a}
\eeq
\beq
    a_\mu^{(k)}(x_1,\ldots,x_{k-1})
        =\begin{vmatrix}
      x_1^{\mu_1} & x_1^{\mu_2} & \ldots & x_1^{\mu_n} \\
             \hdotsfor{4} \\
      x_{k-1}^{\mu_1} & x_{k-1}^{\mu_2} & \ldots & x_{k-1}^{\mu_n} \\
      \mu_1^{n-k} & \mu_2^{n-k} & \ldots & \mu_n^{n-k} \\
             \hdotsfor{4} \\
            \mu_1 & \mu_2 & \ldots & \mu_n \\
            1 & 1 & \ldots & 1
      \end{vmatrix},
\label{eq:amuk}
\eeq
\beq
   a_\d^{(k)}(x_1,\ldots,x_{k-1})=
   \left(\prod_{i=1}^{n-k} i!   \right)
   \left(\prod_{j=1}^{k-1}(x_j-1)^{n-k+1}\right)
   \left( \prod_{i<j}(x_i-x_j)\right).
\label{eq:adk}
\eeq
\end{subequations}
\end{lemma}

{\bf Proof} is given by induction in $(n-k)$. For $k=n$ the statement is
obvious. Then, set $x_{k-1}=1+\varepsilon$ in \Ref{eq:amuk} and expand
$(1+\varepsilon)^{\mu_j}$
in the $k^\text{th}$ row in powers of
$\varepsilon$:
\be
   (1+\varepsilon)^\mu=1+\binom{\mu}{1}\varepsilon+\cdots
   +\binom{\mu}{n-k+1}\varepsilon^{n-k+1}+\cdots
\ee

The terms containing $\mu$ in powers $\leq(n-k)$
are cancelled from the determinant by
subtracting from $k^\text{th}$ row a linear combination of the lower
rows. As a result, we get $a_\mu^{(k)}\sim
\varepsilon^{n-k+1}a_\mu^{k-1}/(n-k+1)!$.
The factor $\varepsilon^{n-k+1}/(n-k+1)!$ is cancelled then by the same factor
in the expansion of the expression \Ref{eq:adk} for $a_\d^{(k)}$.
\endproof

In particular, for $k=1$ we get from \Ref{eq:restr_schur} the formula
for the normalization factor
\beq
    s_\la(\one)=\prod_{i<j}\frac{\mu_i-\mu_j}{j-i}=
    \frac{\D_n(\mu)}{\D_n(\d)},
\qquad
    \D_n(\d)=1!2!3!\ldots(n-1)!
\label{eq:snorm}
\eeq

Using \Ref{eq:snorm} we define the normalized Schur functions
as $\bar{s}_\la(\xbf)\equiv s_\la(\xbf)/s_\la(\one)$, so that $\bar{s}_\la(\one)=1$.

\subsection{Operators $H_j$}
Let $\Cbbd[\xbf]^{A_n}$ be the space of {\it antisymmetric} polynomials.
The multiplication by $\D_n(\xbf)$ provides an isomorphism of linear spaces
and the corresponding mapping of the bases
\beq
\D_n(\xbf):\Cbbd[\xbf]^{S_n}\mapsto\Cbbd[\xbf]^{A_n}
:s_\la(\xbf)\mapsto a_\mu(\xbf),
\label{eq:sym-antisym}
\eeq
see \Ref{eq:def_s} and \Ref{eq:vandermonde}.
The polynomials $a_\mu(\xbf)$ given by \Ref{eq:amu} are linear combinations
of monomials $\xbf^{\sigma(\mu)}$
and therefore are eigenfunctions of the following
commuting differential operators in $\Cbbd[\xbf]^{A_n}$:
\beq
\tilde{H}_j=e_j(D_1,\ldots,D_n),\qquad
\tilde{H}_j a_\mu=e_j(\mu)a_\mu.
\label{eq:tildeHops_s}
\eeq
\be
[\tilde{H}_j,\tilde{H}_k]=0,\qquad \forall j,k=1,\ldots,n,
\ee
cf.\ \Ref{eq:Hops_m}.
Conjugating the operators $\tilde{H}_j$ with the factor $\D_n(\xbf)$
we get the family of commuting differential operators $H_j$
in $\Cbbd[\xbf]^{S_n}$ diagonalized by the basis $s_\la(\xbf)$:
\beq
    H_j=\D_n^{-1}(\xbf)e_j(D_1,\ldots,D_n)\D_n(\xbf), \qquad
    H_js_\la=e_j(\mu)s_\la,
\label{eq:Hops_s}
\eeq
\be
[{H}_j,{H}_k]=0,\qquad \forall j,k=1,\ldots,n,
\ee

\subsection{$Q$-operator}
\begin{prop}
The polynomial $q_\la(z)\equiv \bar{s}_\la(z,1,\ldots,1)$
is given by the expression
\beq
    q_\la(z)=\frac{(n-1)!}{(z-1)^{n-1}}\,\phi_\la(z),
\label{eq:qformula_s}
\eeq
where
\beq
    \phi_\la(z)=\sum_{j=1}^n c_j z^{\mu_j}, \qquad
     c_j=\prod_{k\neq j} (\mu_j-\mu_k)^{-1}.
\label{eq:def-phi}
\eeq
\end{prop}

{\bf Proof.}
Using the formulae \Ref{eq:s=a/a} and \Ref{eq:snorm} we get
\beq
    q_\la(z)= \bar{s}_\la(z,1,\ldots,1)
        =\frac{s_\la(z,1,\ldots,1)}{s_\la(1,\ldots,1)}
        =\frac{\D_n(\d)}{\D_n(\mu)}\frac{a_\mu^{(2)}(z)}{a_\d^{(2)}(z)}.
\label{eq:q=ad/ad}
\eeq

Substituting
$a_\d^{(2)}(z)=(z-1)^{n-1}(n-2)!\ldots2!1!$ from \Ref{eq:adk} 
and $\D_n(\d)$ from \Ref{eq:snorm},
and cancelling the factorials
we transform \Ref{eq:q=ad/ad} into
\beq
    q_\la(z)=\frac{(n-1)!}{(z-1)^{n-1}}\frac{a^{(2)}_\mu(z)}{\D_n(\mu)}
         =\frac{(n-1)!}{(z-1)^{n-1}}\,\phi_\la(z), \qquad
        \phi_\la(z)=\frac{a^{(2)}_\mu(z)}{\D_n(\mu)}.
\eeq

Noticing that, by \Ref{eq:amuk},
\beq
   a_\mu^{(2)}(z) = \begin{vmatrix}
             z^{\mu_1} & z^{\mu_2} & \ldots & z^{\mu_n} \\
             \mu_1^{n-2} & \mu_2^{n-2} & \ldots & \mu_n^{n-2} \\
             \hdotsfor{4} \\
            \mu_1 & \mu_2 & \ldots & \mu_n \\
            1 & 1 & \ldots & 1
          \end{vmatrix}
\eeq
we expand the above determinant along the first row and,
cancelling the arising Vandermonde determinants, obtain
\Ref{eq:def-phi}.
\endproof

The polynomial $\phi_\la(z)$ defined by \Ref{eq:def-phi}
is a linear combination of the monomials
$z^{\mu_j}$, $j=1,\ldots,n$, and, therefore, satisfies the differential equation
\beq
    \prod_{j=1}^n (z\dd_{z}-\mu_j)\phi=0.
\label{eq:diffeq_phi}
\eeq

The polynomiality of $q_\la(z)$ implies, by virtue of \Ref{eq:qformula_s},
that the polynomial $\phi_\la(z)$ is divisible by $(z-1)^{n-1}$.
To verify this property directly note that the divisibility is
equivalent to the relations
\beq
   \sum_{j=1}^n \mu_j^k c_j=0, \qquad k=0,\ldots,n-2,
\label{eq:eqs_cj1}
\eeq
which, taken together with the normalization condition
$q_\la(1)=1$ expressed as
\beq
    \sum_{j=1}^n \mu_j^{n-1} c_j=1,
\label{eq:eqs_cj2}
\eeq
follow from the identities
\be
    \sum_{j=1}^n\frac{\mu_j^l}{\prod_{k\neq j}(\mu_j-\mu_k)}=\delta_{l,n-1},
    \quad l=0,\ldots,n-1,
\ee
which are obtained by setting $t=0$ in the partial fraction decomposition
\be
    \frac{t^{l+1}}{\prod_{k=1}^n(t-\mu_k)}
    =\delta_{l,n-1}+\sum_{j=1}^n
\frac{1}{t-\mu_j}\,\frac{\mu_j^{l+1}}{\prod_{k\neq j}(\mu_j-\mu_k)}.
\ee

In fact, the above properties can be used to characterize $\phi_\la(z)$.

\begin{prop}
The solution $\phi_\la(z)$ to \Ref{eq:diffeq_phi}
is uniquely characterized by the conditions $[\dd_z^k\phi]_{z=1}=0$, $k=0,\ldots,n-2$,
which ensure that the polynomial $\phi_\la(z)$ is divisible by $(z-1)^{n-1}$, and
by $[\dd_z^{n-1}\phi]_{z=1}=1$, which is equivalent to the normalization
condition $q_\la(1)=1$.
\end{prop}
{\bf Proof.} The general solution to \Ref{eq:diffeq_phi} is
$\phi=\sum_{j=1}^n c_jz^{\mu_j}$. The divisibility by $(z-1)^{n-1}$
condition is equivalent then to the equations \Ref{eq:eqs_cj1},
and the normalization condition to the equation \Ref{eq:eqs_cj2}.
Solving the system of linear equations \Ref{eq:eqs_cj1}--\Ref{eq:eqs_cj2}
for $c_j$ via Cramer's formula we recover the unique solution \Ref{eq:def-phi}.
\endproof

Note that the differential equation \Ref{eq:diffeq_phi} can be written also as
\be
   \left[(z\dd_z)^n+\sum_{k=1}^n (-1)^k h_k(\la) (z\dd_z)^{n-k}\right]\phi=0
\ee
where $h_k(\la)=e_k(\mu)$ are the eigenvalues of the operators $H_k$
\Ref{eq:Hops_s}.
The corresponding differential equation for $q_\la(z)$ is
\beq
   \left[Z^n+\sum_{k=1}^n (-1)^k h_k(\la) Z^{n-k}\right]q_\la(z)=0,
\qquad Z=z\left(\dd_z+\frac{n-1}{z-1}\right).
\label{eq:deffeq_q_s}
\eeq

Note that $q_\la(z)$ is the only, up to a coefficient, polynomial
solution to \Ref{eq:deffeq_q_s}.

As usual, we define the operator $Q_z$ in $\Cbbd[\xbf]^{S_n}$ through
its eigenvectors $\bar{s}_\la(\xbf)$ and eigenvalues
$q_\la(z)$.

\begin{theo}\label{th:Q_s}
Let $z>1$ and $y_1<y_2<\ldots<y_n$. Then for any $f\in\Cbbd[\xbf]^{S_n}$
the value $[Q_zf](\ybf)$ is given by the integral
\beq
    [Q_zf](\ybf)=(n-1)!\frac{(z-1)^{n-1}}{\D_n(\ybf)}
    \int_{\Omega_x} d\xbf\, \d(x_1\ldots x_n-zy_1\ldots y_n)
    \D_n(\xbf)f(\xbf),
\label{eq:Qcal}
\eeq
where the integration domain $\Omega_x$ is defined by the inequalities
\beq
    0<y_1<x_1<y_2<\ldots<x_{n-1}<y_n<x_n.
\label{eq:def_Omega}
\eeq
\end{theo}

The above formula is obtained by setting $g=1$ in 
the formula (4.2) in \cite{KMS}.
Here we present an elementary and independent proof.

{\bf Proof.} From \Ref{eq:def_Omega} it follows that $x_j>y_j$ for $j=1,\ldots,n$,
hence $x_1\ldots x_n>y_1\ldots y_n$. The condition $z>1$ ensures then that
the support of the delta function in \Ref{eq:Qcal} has a non-empty intersection with
the domain $\Omega_x$.

Since the Schur functions $s_\la(\xbf)$ form a basis in
$\Cbbd[\xbf]^{S_n}$ it is sufficient to verify \Ref{eq:Qcal} for $f=s_\la$:
\beq
    (n-1)!\frac{(z-1)^{n-1}}{\D_n(\ybf)}
    \int_{\Omega_x} d\xbf\, \d(x_1\ldots x_n-zy_1\ldots y_n)
    \D_n(\xbf)s_\la(\xbf)=q_\la(z)s_\la(\ybf).
\eeq

Using the correspondence \Ref{eq:sym-antisym} of the symmetric and antisymmetric polynomials
as well as the formula \Ref{eq:qformula_s}
we reduce the task to verifying the equality
\beq
    \int_{\Omega_x} d\xbf\, \d(x_1\ldots x_n-zy_1\ldots y_n)a_\mu(\xbf)
    =a_\mu(\ybf)\phi_\la(z).
\label{eq:int_amu}
\eeq

Expanding the determinantal expression \Ref{eq:amu} for $a_\mu(\xbf)$
along the last row we get
\beq
    a_\mu(\xbf)=\sum_{k=1}^n (-1)^{k+n}x_n^{\mu_k}\det M^{(k)}
\label{eq:expand_amu}
\eeq
where the matrix $M^{(k)}$ is
\beq
    [M^{(k)}]_{ij}=x_i^{\mu_k}, \qquad i=1,\ldots,n-1; \qquad
     j=1,\ldots,\hat{k},\ldots,n
\label{eq:matrixMindices}
\eeq

Integrating \Ref{eq:expand_amu} in the variable $x_n$ with the
delta-function factor
$\d(x_1\ldots x_n-zy_1\ldots y_n)$ replaces $x_n^{\mu_k}$ with
the factor
\be
   \frac{(zy_1\ldots y_n)^{\mu_k}}{(x_1\ldots x_{n-1})^{\mu_k+1}},
\ee
hence
\beq
    \int dx_n\,\d(\ldots)\,a_\mu(\xbf)=
     \sum_{k=1}^n (-1)^{k+n}(zy_1\ldots y_n)^{\mu_k}\,
     x_n^{\mu_k}\det\{x_i^{\mu_j-\mu_k-1}  \}
\eeq
where the matrix indices $i$,$j$ run like in \Ref{eq:matrixMindices}.
Further integration is performed independently for each row of the determinant:
\beq
    \int_{y_i}^{y_{i+1}}dx_i\,x_i^{\mu_{jk}-1}=
    \frac{y_{i+1}^{\mu_{jk}}-y_i^{\mu_{jk}}}{\mu_j-\mu_k}
\eeq
(the difference $\mu_{jk}\equiv\mu_j-\mu_k$ is never 0, so there
are no logarithms). According to \Ref{eq:def-phi}
the product of the factors $(\mu_j-\mu_k)^{-1}$
produces the coefficient $(-1)^{n-1}c_k$, and the left-hand side of
\Ref{eq:int_amu} is transformed into
\beq
    \sum_{k=1}^n (-1)^{k-1}(zy_1\ldots y_n)^{\mu_k}c_k
    \det\{y_{i+1}^{\mu_{jk}}-y_i^{\mu_{jk}} \}.
\label{eq:intamu_intermediate}
\eeq

The determinant of order $(n-1)$ in \Ref{eq:intamu_intermediate}
is transformed into a determinant of order $n$ in the following way.
Let $t_{i}^{j}=y_i^{\mu_{jk}}$. Then
\begin{multline}
    \begin{vmatrix}
      t_{2}^{1}-t_{1}^{1} & \ldots & t_{2}^{k-1}-t_{1}^{k-1} & t_{2}^{k+1}-t_{1}^{k+1} & \ldots & t_{2}^{n}-t_{1}^{n} \\
      t_{3}^{1}-t_{2}^{1} & \ldots & t_{3}^{k-1}-t_{2}^{k-1} & t_{3}^{k+1}-t_{2}^{k+1} & \ldots & t_{3}^{n}-t_{2}^{n} \\
             \hdotsfor{6} \\
      t_{n}^{1}-t_{n-1}^{1} & \ldots & t_{n}^{k-1}-t_{n-1}^{k-1} & t_{n}^{k+1}-t_{n-1}^{k+1} & \ldots & t_{n}^{n}-t_{n-1}^{n}
    \end{vmatrix} \\
    = (-1)^{k-1}
    \begin{vmatrix}
       t_{1}^{1} & \ldots & t_{1}^{k-1} & 1 & t_{1}^{k+1} & \ldots & t_{1}^{n} \\
       t_{2}^{1} & \ldots & t_{2}^{k-1} & 1 & t_{2}^{k+1} & \ldots & t_{2}^{n} \\
             \hdotsfor{7} \\
       t_{n}^{1} & \ldots & t_{n}^{k-1} & 1 & t_{n}^{k+1} & \ldots & t_{n}^{n}
    \end{vmatrix}.
\label{eq:matrix_id}
\end{multline}

To prove the matrix identity \Ref{eq:matrix_id} take its right-hand side
and replace the $i^\text{th}$ row,
for $i$ from $2$ to $n$, with the difference of the $i^\text{th}$ and $(i-1)^\text{th}$
row, so that $k^\text{th}$ column becomes all zeroes except the element $(1k)$.
Expanding the determinant in the $k^\text{th}$ column produces the left-hand side of
\Ref{eq:matrix_id} which coincides with the determinant in \Ref{eq:intamu_intermediate}.

Using the factors $(zy_1\ldots y_n)^{\mu_k}$ in \Ref{eq:intamu_intermediate}
we transform the right-hand side of \Ref{eq:matrix_id} to the familiar form
\Ref{eq:amu} for $a_\mu(\ybf)$ and finally reproduce the right-hand side of
\Ref{eq:int_amu}.
\endproof

\subsection{Separation of variables}

\begin{theo}
For the separating operator $\mathcal{S}_n$ defined
by the formula \Ref{eq:def_S}
and satisfying \Ref{eq:actSP}
for $\bar{P}_\la\equiv \bar{s}_\la$
there exist uniquely defined operators $\Acal_k$ satisfying
the relations \Ref{eq:rhoQ=Arho}. The value of
$[\Acal_{k}f](\tilde{\ybf})$ is given by the integral
\begin{multline}
    [\Acal_{k}f](\tilde{\ybf})
    =\frac{(-1)^{k-1}(n-1)!}{(n-k)!%
(z_{k}-1)^{n-1}\D_{k-1}(\tilde{\ybf})\prod_{j=1}^{k-1}(\tilde{y}_j-1)^{n-k+1}} \\
    \times\int_{\tilde{\Omega}_x} d\tilde{\xbf}\,
    \d(\tilde{x}_1\ldots \tilde{x}_{k}-z_{k}\tilde{y}_1\ldots \tilde{y}_{k-1})\,
    \D_{k}(\tilde{\xbf})\,
    \prod_{j=1}^k (\tilde{x}_j-1)^{n-k}
    f(\tilde{\xbf})
\label{eq:Acal_s}
\end{multline}
where $\tilde{\xbf}=(\tilde{x}_1,\ldots,\tilde{x}_{k})$,
$\tilde{\ybf}=(\tilde{y}_1,\ldots,\tilde{y}_{k-1})$,
and the integration domain is
\be
\tilde{\Omega}_x:\quad 1<\tilde{x}_1<\tilde{y}_1<\ldots<\tilde{y}_{k-1}<\tilde{x}_k\,.
\ee
\end{theo}

{\bf Proof.} The formula \Ref{eq:Acal_s} can be obtained directly from
the formula (6.18) from \cite{KMS} by setting $g=1$. The proof
given below is basically a simplified proof of the Proposition 6.1 from 
\cite{KMS}. The mismatch of signs with respect to \cite{KMS}
is due to a changed definition of $\D_n$.

Let us evaluate $\rho_{k-1}Q_{z_k}$ in the left-hand side
of \Ref{eq:rhoQ=Arho} using the definition \Ref{eq:def_rhok}
of $\rho_k$ and the integral formula \Ref{eq:Qcal} for $Q_z$.
The operator $\rho_{k-1}$ sets $(n-k+1)$ of the variables $y_j$
in \Ref{eq:Qcal} to the unit values. Since $[Q_{z_k}f](\ybf)$
is a symmetric polynomial it does not matter which of $y_j$ we choose
to fix. Let us set
\begin{alignat*}{2}
    y_j&=1+\varepsilon v_j, &\qquad j&=1,\ldots,n-k+1, \\
    y_{n-k+j+1}&=\tilde{y}_j, &\qquad j&=1,\ldots, k-1, \\
    x_j&=1+\varepsilon u_j, &\qquad j&=1,\ldots,n-k, \\
    x_{n-k+j}&=\tilde{x}_j, &\qquad j&=1,\ldots,k
\end{alignat*}
and take the limit $\varepsilon\rightarrow0$.
Due to the inequalities \Ref{eq:def_Omega} the variables
$x_1,\ldots,x_{n-k}$ are pinched between $y$'s and, therefore,
forced to tend to 1 as well, which accounts for the factor $\rho_{k+1}$
in the right-hand side $\Acal_k\rho_k$
of \Ref{eq:rhoQ=Arho}. To calculate the kernel
of the integral operator $\Acal_k$ we observe that in the limit
$\varepsilon\rightarrow0$
\begin{align*}
  \D_n(\xbf)\sim\varepsilon^{\frac{(n-k)(n-k-1)}{2}}
     \D_{n-k}(\ubf)\D_k(\tilde{\xbf})
     \prod_{j=1}^k(1-\tilde{x}_j)^{n-k}, \\
  \D_n(\ybf)\sim\varepsilon^{\frac{(n-k+1)(n-k)}{2}}
     \D_{n-k+1}(\vbf)\D_{k-1}(\tilde{\vbf})
     \prod_{j=1}^{k-1}(1-\tilde{y}_j)^{n-k+1}.
\end{align*}

Since $dx_1\ldots dx_{n-k}\sim\varepsilon^{n-k} du_1\ldots du_{n-k}$,
the factors $\varepsilon$ cancel completely with those from 
$\D_n(\xbf)$ and $\D_n(\ybf)$.
The integration in $u_j$ produces the constant factor
\beq
    \int_{v_1}^{v_2}du_1\cdots \int_{v_{n-k}}^{v_{n-k+1}} du_{n-k}\,
    \D_{n-k}(\ubf)=
    \frac{(-1)^{n-k}}{(n-k)!}\D_{n-k+1}(\vbf).
\label{eq:int_uv}
\eeq

The formula \Ref{eq:int_uv} is proved by integrating independently in $u_i$
the $i$-th row of the determinant representing $\D_{n-k}(\ubf)$ and using
an identity similar to \Ref{eq:matrix_id}. Collecting then all the factors and
coefficients we finally get \Ref{eq:Acal_s}.
\endproof

A peculiar feature of the Schur functions ($\a=1$) case is a beautiful
inversion formula for the operator $\Scal_n$.

\begin{theo}
The inverse of $\Scal_n$ is the differential operator on
$\Cbbd[\xbf]^{S_n}$ given by the formula
\beq
    \Scal_n^{-1}=\frac{(-1)^{\frac{n(n-1)}{2}}}{[(n-1)!]^n}
    \frac{\D_n(\d)}{\D_n(\xbf)}
     \circ K_n\circ
     \prod_{k=1}^n(x_k-1)^{n-1},
\label{eq:Sinv_schur}
\eeq
where $K_n$ is the differential operator
\beq
     K_n=\det\{D_i^{\d_j}\}
     =\prod_{i<j}(D_i-D_j),
     \qquad D_i\equiv x_i\frac{\dd}{\dd x_i}.
\label{eq:def_K}
\eeq
\end{theo}

Since $\Scal_n^{-1}$ is a differential operator it is convenient to
identify $\xbf$ and $\zbf$ variables and assume that $\Scal_n^{-1}$
acts in $\Cbbd[\xbf]^{S_n}$.
The formula for $\Scal_n^{-1}$ was mentioned in \cite{KS4} without a proof
(and in a slightly different notation).
Here we present a detailed derivation.

{\bf Proof.}
It is sufficient to verify that
\beq
   S_n^{-1}:\quad \prod_{j=1}^n q_\la(x_j)\mapsto
   \bar{s}_\la(\xbf).
\eeq

Using the formulae \Ref{eq:def_s} for $s_\la(\xbf)$,
\Ref{eq:vandermonde} for $\D_n(\xbf)$,
\Ref{eq:snorm} for $s_\la(\one)$,
\Ref{eq:qformula_s} for $q_\la(z)$,
and \Ref{eq:def-phi} for $\phi_\la(z)$
we reduce the task to proving the equality
\beq
    K_n: \prod_{i=1}^n \phi_\la(x_i)
    \mapsto (-1)^{\frac{n(n-1)}{2}} \frac{a_\mu(\xbf)}{\D_n(\mu)}.
\label{eq:Kphi}
\eeq

Using the determinantal representation \Ref{eq:def_K}
for the operator $K_n$ we represent $K_n\prod\phi_\la(x_i)$
as the determinant of the matrix $M$:
\beq
    M_i^j=D_i^{\d_j}\phi_\la(x_i)
         =D_i^{\d_j}\sum_{k=1}^n c_k x_i^{\mu_k}
         =\sum_{k=1}^n x_i^{\mu_k}\cdot c_k\mu_k^{\d_j}\,.
\eeq

The resulting expression is recognized as a product of three matrices,
therefore
\begin{align}
    K_n \prod_{i=1}^n \phi_\la(x_i) &= \det[M_i^j]
    =\det[x_i^{\mu_j}]\cdot \det[c_i\d_{ij}]\cdot \det[\mu_i^{\d_j}] \\
    &=a_\mu(\xbf)\cdot(c_1\cdots c_n)\cdot \D_n(\mu). \notag
\end{align}
To obtain \Ref{eq:Kphi} and thus conclude the proof it remains
to substitute into the above formula the values of $c_n$ from \Ref{eq:def-phi}.
\endproof

The differential operator $K_n$ appeared also in \cite{BG62}
though in a different context.

\subsection{Lifting operator}

\begin{prop}
The operator $Q_z$ taken at $z=0$ admits a unique factorization
\Ref{eq:factQ0} into the projector $\Pcal$ \Ref{eq:Pcal}
and the operator $Q_0^\prime$ \Ref{eq:Q0prime}:
\beq
Q_0^\prime: \bar{s}_{\la_1\ldots\la_{n-1}}(x_1,\ldots,x_{n-1})
\mapsto \bar{s}_{\la_1\ldots\la_{n-1}0}(x_1,\ldots,x_{n}).
\label{eq:Q0prime_s_lift}
\eeq

For $y_1<y_2<\cdots<y_n$ the value $[Q_0^\prime f](\ybf)$
is given by the integral
\beq
  [Q_0^\prime f](\ybf)=(-1)^{n-1}\frac{(n-1)!}{\D_n(\ybf)}
  \int_{\Omega_x^\prime} d\xbf'\,
  \D_{n-1}(\xbf')f(\xbf')
\label{eq:Qp0prime_s}
\eeq
where $\xbf'=(x_1,\ldots,x_{n-1})$ and
\beq
  \Omega_x^\prime:\quad 0<y_1<x_1<y_2<\ldots<x_{n-1}<y_n.
\eeq
\end{prop}

{\bf Proof.} We cannot set $z=0$ directly in the integral
\Ref{eq:Qcal} because of the restriction $z>1$.
Instead, we shall use the definition of $Q_z$ in terms of its eigenvectors
$s_\la(\xbf)$ and eigenvalues $q_\la(z)$.

In what follows we have to be careful to distinguish
the partitions of length $n$ and $n-1$.
Let $\la=(\la_1,\ldots,\la_n)$ and $\la'=(\la_1,\ldots,\la_{n-1})$.
Respectively,
$\mu=(\mu_1,\ldots,\mu_n)$ and $\mu'=(\mu_1,\ldots,\mu_{n-1})$,
where $\mu_i=\la+n-i$, $\mu_i^\prime=\la_i+n-1-i$. Note that
$\mu_i=\mu_i^\prime+1$, for $i=1,\ldots,n-1$.

{}From \Ref{eq:qformula_s} and \Ref{eq:def-phi} it follows that
$q_\la(0)$ is nonzero only for $\mu_n\equiv\la_n=0$. In this case
\beq
    q_{\la_1\ldots\la_{n-1}0}(0)=
    \frac{(n-1)!}{\prod_{i=1}^{n-1}\mu_i}.
\label{eq:q0schur}
\eeq

Similarly, from \Ref{eq:def_s} and \Ref{eq:amu} we conclude that
$\Pcal s_\la$ does not vanish only for $\la_n=0$, and
\beq
    s_{\la_1\ldots\la_{n-1}0}(x_1,\ldots,x_{n-1},0)
    =s_{\la_1\ldots\la_{n-1}}(x_1,\ldots,x_{n-1}).
\eeq

Using \Ref{eq:snorm}
we get, respectively, for the normalized Schur functions
\beq
    \bar{s}_{\la_1\ldots\la_{n-1}0}(x_1,\ldots,x_{n-1},0)
    =\frac{(n-1)!}{\prod_{i=1}^{n-1}\mu_i}
    \bar{s}_{\la_1\ldots\la_{n-1}}(x_1,\ldots,x_{n-1}).
\eeq

{}From the last formula together with \Ref{eq:q0schur}
the factorization \Ref{eq:factQ0} and the formula
\Ref{eq:Q0prime_s_lift} follow immediately.

The integral formula \Ref{eq:Qp0prime_s} can be obtained, in principle,
by setting $g=1$ in the formula (7.10) from \cite{KMS} and taking into account
the different definition of $\D_n$. We provide, however, an independent proof.
It is sufficient to verify \Ref{eq:Qp0prime_s} on the basis
$f=\bar{s}_{\la'}(\xbf')$. Using \Ref{eq:def_s} and \Ref{eq:snorm},
and cancelling the arising coefficients we reduce the problem to verifying the
identity
\beq
    \int_{\Omega_x^\prime} d\xbf'\, a_{\mu^\prime}(\xbf')
    =\frac{(-1)^{n-1}}{\prod_{i=1}^{n-1}\mu_i}\,
    a_\mu(\ybf).
\eeq

To prove the last identity we integrate independently
in $x_i$ from $y_i$ to $y_{i+1}$ each row of the determinant representing
$a_{\mu^\prime}(\xbf')$, and then use $\mu_i=\mu_i^\prime+1$ and
a variant of the determinantal identity \Ref{eq:matrix_id}.
\endproof

\section{Concluding remarks}
\noindent
We have considered three important standard bases in the linear space of
symmetric polynomials in $n$ variables. These bases have been related to
the special cases of the Jack polynomials which, in their turn, solve
the famous quantum integrable system (Calogero-Sutherland model).
The main objective was to demonstrate the main features of the
(quantum) separation of variables, which is designed to
produce explicit factorization and representation for the
multivariate special functions through the application
of suitable (integral) operators.

The method of quantum separation of variable has its
counterpart in the classical Hamiltonian mechanics which is
applicable to a wide class of Liouville integrable
systems. Usually, it is the classical system that gets
separated in the first place, followed by the problem
of quantization. For the three bases considered in this paper
such approach is valid only for the $E_\la$ one because
the Jack's parameter $\alpha$ is proportional to the Plank
constant. The two other bases, $m_\la$ and $s_\la$, do not
have classical analogs, so that the separation of variables
that we have found for them in this paper is purely quantum.

In this paper we considered polynomials symmetric under permutation group
$S_n$, or Weyl group of the root system $A_{n-1}$.
It is a challenging problem to find the factorizing operators
for the bases of polynomials symmetric with respect to Weyl groups
corresponding to other root systems, such as, for example,
$B_n$, $C_n$, or $D_n$ \cite{Weyl39,Lit40}.



\begin{thebibliography}{12}

\bibitem{KMS}
V.B.~Kuznetsov, V.V.~Mangazeev and E.K.~Sklyanin,
{\em $Q$-operator and factorized separation
chain for Jack polynomials}.
\newblock Indag. Mathem.~{\em 14}(3,4) (2003), pp.\ 451--482.
ArXiv: http://arxiv.org/abs/math.CA/0306242

\bibitem{Macd}
I.G.~Macdonald,
{\it Symmetric functions and Hall polynomials},
Oxford Univ. Press, Oxford, Second Edition (1995).

\bibitem{KS4} V.B.~Kuznetsov and E.K.~Sklyanin,
{\it Separation of variables in the
$A_2$ type Jack polynomials,} RIMS
Kokyuroku, {\em 919} (1995), 27--34.
ArXiv: http://arxiv.org/abs/solv-int/9508002

\bibitem{BG62} F.A.~Berezin and I.M.~Gelfand,
{\em Some remarks on the theory of spherical functions on symmetric Riemannian manifolds},
Amer. Math. Soc. Transl. (2) {\bf 21} (1962) 193--238.

\bibitem{Weyl39} H.~Weyl, {\it The Classical Groups, Their Invariants and
Representations}, Princeton University Press, Princeton (1939).

\bibitem{Lit40} D.E.~Littlewood,
{\em The theory of group characters and matrix representations of
groups}, Oxford Univ.\ Press, Oxford (1940).
\end{thebibliography}
\end{document}